\documentclass[a4paper,12pt,reqno]{amsart}

\usepackage{amsmath}
\usepackage{tikz-cd}
\usepackage{amscd}
\usepackage{amssymb}
\usepackage{mathrsfs}
\usepackage[left=2.5cm,right=2.5cm,bottom=3cm,top=3cm]{geometry}
\newtheorem{thm}{Theorem}[section]

\newtheorem{prop}[thm]{Proposition}
\theoremstyle{definition}

\newtheorem{rem}[thm]{Remark}

\numberwithin{equation}{section}
\newcommand{\dd}{\sqrt{-1}\partial\bar{\partial}}

\begin{document}

\title[Compact K\"ahler manifolds with a type-III solution to KRF]{Compact K\"ahler manifolds with a type-III solution to the K\"ahler-Ricci flow}

\author{Yashan Zhang}

\thanks{Partially supported by National Natural Science Foundation of China (No. 12371057) and Natural Science Foundation of Hunan Province (No. 2024JJ2006)}

\address{Hunan University, Changsha, China}
\email{yashanzh@hnu.edu.cn}

\begin{abstract}
On a compact K\"ahler manifold of a semi-ample canonical line bundle, Tosatti and Y.G. Zhang \cite{ToZyg} developed an almost complete classification on infinite-time singularity type of the K\"ahler-Ricci flow, which is complete in surface case. In this paper, we present a complete classification for the projective case. Precisely, we give a characterization for a projective manifold of a semi-ample canonical line bundle admitting a type-III solution to the K\"ahler-Ricci flow. Partial extensions to the compact K\"ahler case are also obtained, which in particular contain a complete classification in the case of compact K\"ahler threefolds. These results may be regarded as extensions of previous structure theorems for a compact K\"ahler manifold of nonpositive holomorphic curvatures.
\end{abstract}

\maketitle

\section{Introduction}
We continue our study in \cite{Zys20a,Zys20b} on the singularity type of long-time solutions to the K\"ahler-Ricci flow on compact K\"ahler manifolds. Our main interest is the \emph{relation between infinite-time singularity type (more precisely, long time behavior of curvature) of the K\"ahler-Ricci flow and the geometric/complex structure of the underlying compact K\"ahler manifold}. 
\par Thanks to the maximal existence time theorem of the K\"ahler-Ricci flow \cite{C,TZo,Ts}, the existence of long time solutions is equivalent to nefness of the canonical line bundle. In this paper, we will further restrict our discussions on compact K\"ahler manifold of a semi-ample canonical line bundle.
\par Let $X$ be an $n$-dimensional compact connected K\"ahler manifold of a semi-ample canonical line bundle $K_X$, and so the Kodaira dimension $kod(X)\in\{0,1,\ldots,n\}$. Given an arbitrary K\"ahler metric $\omega_0$ on $X$, let $\omega=\omega(t)_{t\in[0,\infty)}$ be the solution to the (normalized) K\"ahler-Ricci flow
\begin{equation}\label{nkrf}
\partial_t\omega(t)=-Ric(\omega(t))-\omega(t)
\end{equation}
running from $\omega_0$. Recall from \cite{Ha93} that we say a long time solution $\omega(t)$ to the K\"ahler-Ricci flow \eqref{nkrf} is of type \textrm{III} if
\begin{equation}\label{defn_ii}
\limsup_{t\to\infty}\left(\sup_{X}|Rm(\omega(t))|_{\omega(t)}\right)<\infty\nonumber
\end{equation}
and of type  \textrm{IIb} if
\begin{equation}\label{defn_iii}
\limsup_{t\to\infty}\left(\sup_{X}|Rm(\omega(t))|_{\omega(t)}\right)=\infty\nonumber.
\end{equation}
For $X$ with semi-ample canonical line bundle, we let
\begin{equation}\label{eq-semiample}
f:X\to X_{can}\subset \mathbb{CP}^N
\end{equation}
be the semi-ample fibration with connected fibers induced by pluricanonical system of $K_X$, where $X_{can}$, a $kod(X)$-dimensional irreducible normal projective variety, is the canonical model of $X$, and $V\subset X_{can}$ be the singular set of $X_{can}$ together with the critical values of $f$, and $X_y=f^{-1}(y)$ be a smooth fiber for $y\in X_{can}\setminus V$. An almost complete classification of the singularity type of the K\"ahler-Ricci flow on such manifold $X$ has been developed by Tosatti-Y.G. Zhang \cite[Theorems 1.5 and 1.6]{ToZyg} (also see \cite{FZ,Gill,HeTo,SW} for some early works), which leaves open the following case:\\

\begin{itemize}
\item[($\star$)]\,\,\,\,\emph{$0<kod(X)<n$ ($n\ge3$), the generic fiber $X_y$ is a finite \'etale quotient of a complex torus, and $V\neq\emptyset$}.\\
\end{itemize}

\par In general, as expected by Tosatti-Y.G. Zhang \cite[Section 1, page 2930]{ToZyg} and confirmed by our previous work \cite[Theorem 1.4]{Zys20a}, the singularity type of the K\"ahler-Ricci flow in case ($\star$) (and hence in the semi-ample canonical line bundle case) does not depend on the choice of the initial metric (also see \cite[Theorem 1.6]{FZys} for a stronger local version). Actually, it is further conjectured by Tosatti \cite[Conjecture 6.7]{To.kawa} that the singularity type of the K\"ahler-Ricci flow on a compact K\"ahler manifold of a nef canonical line bundle does not depend on the choice of the initial metric, which is affirmatively solved by Wondo and Z. Zhang \cite{WoZz} very recently. Given these results, one sees that the infinite-time singularity type of the K\"ahler-Ricci flow depends only on the complex structure of the underlying compact K\"ahler manifold, and hence it is a very natural problem to classify the structure of $X$ in terms of infinite-time singularity type of the K\"ahler-Ricci flow, as re-emphasized in \cite[Subsection 6.3]{To.kawa} and \cite[Subsection 8.2]{To24}.

On the other hand, as will be explained in Appendix A, a compact K\"ahler manifold of nonpositive bisectional curvature (or a projective manifold of a semi-ample canonical line bundel and nonpositive holomorphic sectional curvature) admits a type-III solution to the K\"ahler-Ricci flow, and hence \emph{admitting a type-III solution to the K\"ahler-Ricci flow may be regarded as a natural curvature condition (strictly) generalizing nonpositive bisectional curvature or nonpositive holomorphic sectional curvature}. Geometric properties of nonpositive holomorphic curvatures have been widely studied in the past decades, see \cite{HLW16,HLWZ,Liu,WuZh,Zhe95} and references therein.\\

The first main result of this paper is a complete classification of singularity types in the projective case.

\begin{thm}\label{thm-proj-a}
Let $X$ be an $n$-dimensional connected projective manifold of a semi-ample canonical line bundle $K_X$ and $1\le kod(X)\le n-1$. The followings are equivalent.
\begin{itemize}
\item[(1)] $X$ admits a type-III solution to the K\"ahler-Ricci flow. 
\item[(2)] There exists a finite \'etale cover $Z\to X$ such that $Z$ is a connected projective manifold of a semi-ample canonical line bundle $K_Z$, $kod(Z)=kod(X)$, and the semi-ample fibration of $Z$ is a holomorphic submersion over a connected projective manifold with every fiber being a finite \'etale quotient of a complex torus.
\end{itemize}
\end{thm}

This settles the classification of the singularity type of the K\"ahler-Ricci flow in case ($\star$) in the \emph{projective} case. Moreover, up to a further finite \'etale cover of $Z$, we may require that the semi-ample fibration $h$ is a holomorphic submersion with complex torus fibers, see Remark \ref{rem-torusfibration}. \\

It is natural to extend the above conclusion in projective case to the general compact K\"ahler case. However, our proof for Theorem \ref{thm-proj} will make essential use of structure properties of \emph{projective} manifolds of a relatively nef cotangent bundle (see \cite{DPS,Ho}), whose K\"ahler analog seems currently unclear. 

Instead, for the K\"ahler case we will restrict to the particular setting of one (numerical) Kodaira dimensional and use a different method to prove the following structure result, which is our second main result.

\begin{thm}\label{thm-kahler-a}
Let $X$ be an $n$-dimensional compact connected K\"ahler manifold of a nef canonical line bundle $K_X$ and the numerical Kodaira dimension $\nu(X)=1$. The followings are equivalent.
\begin{itemize}
\item[(1)] $X$ admits a type-III solution to the K\"ahler-Ricci flow. 
\item[(2)] $K_X$ is semi-ample, and there exists a finite \'etale cover $Z\to X$ such that $Z$ is a compact connected K\"ahler manifold of a semi-ample canonical line bundle $K_Z$, $kod(Z)=1$, and the semi-ample fibration of $Z$ is a holomorphic submersion over a connected closed Riemann surface with every fiber being a complex torus.
\end{itemize}
\end{thm}

We note that in the above one Kodaira dimension case, the base spaces $X_{can}$ is a closed Riemann surface and there are only finitely many critical values of $f$ on $X_{can}$. In this special case, assuming the existence of a type-III solution, it is possible to give a more explicit description for the singular fibers of $f:X\to X_{can}$, which can be applied to recover a necessary condition for the existence of type-III solution proved in our previous work \cite[Theorem 1.3]{Zys20b}, see Remark \ref{rem-onekod} for a full discussion.\\

In higher Kodaira dimensions, we shall settle the case of K\"ahler threefold of Kodaira dimension two.
\begin{thm}\label{thm-kahler-three-a}
Let $X$ be a $3$-dimensional compact connected K\"ahler manifold of a semi-ample canonical line bundle $K_X$ and Kodaira dimension $kod(X)=2$. The followings are equivalent.
\begin{itemize}
\item[(1)] $X$ admits a type-III solution to the K\"ahler-Ricci flow. 
\item[(2)] There exists a finite \'etale cover $Z\to X$ such that $Z$ is a connected compact K\"ahler threefold of a semi-ample canonical line bundle $K_Z$, $kod(Z)=2$, and the semi-ample fibration of $Z$ is a holomorphic submersion over a smooth connected projective surface with every fiber being a complex torus.
\end{itemize}
\end{thm}

\begin{rem}
Combining our Theorems \ref{thm-kahler-a} and \ref{thm-kahler-three-a} and Tosatti-Y.G. Zhang's work \cite[Theorem 1.5, $\kappa=0,\, n$]{ToZyg}, and abundence theorem for compact K\"ahler threefolds \cite{CHP}, we have a complete classification on singularity type of the long-time K\"ahler-Ricci flow on a compact K\"ahler threefold. 
\end{rem}

\subsection*{Organization} We will first prove several implications of the existence of a type-III solution in the next section, basing on which we will prove Theorems \ref{thm-proj-a}, \ref{thm-kahler-a} and \ref{thm-kahler-three-a} in Sections \ref{sect-proof-proj}, \ref{sect-proof-kahler} and \ref{sect-proof-kahler-three}, respectively. In Appendix A, we will exhibit examples of type-III solution to the K\"{a}hler-Ricci flow under nonpositive holomorphic curvature conditions.

\section{General implications of a type-III solution}
We shall prove two useful general implications of the existence of a type-III solution to the K\"ahler-Ricci flow. 
\subsection{Relative numerical flatness of $T_X^{1,0}$} In this subsection, we assmue $X$ is an $n$-dimensional compact K\"ahler manifold of a semi-ample canonical line bundle $K_X$, and $1\le kod(X)\le n-1$, and $f:X\to X_{can}$ is the semi-ample fibration as in \eqref{eq-semiample}.

\begin{prop}\label{prop-f-flat}
Assume there exists a type-III solution $\omega(t) $ to the K\"ahler-Ricci flow \eqref{nkrf} on $X$. Then $T^{1,0}_X$ is $f$-numerically flat.
\end{prop}
Recall from \cite{DPS,Ho} that a holomorphic vector bundle $E\to X$ over $X$ is nef if $\mathcal O_{\mathbb{P}(E)}(1)$ is nef, and $E$ is numerically flat if both $E$ and $E^*$ are nef. Moreover, $E$ is $f$-nef if $E|_S$ is nef for every reduced irreducible compact analytic subvariety $S\subset X$ with $f(S)=\{pt\}$, and $E$ is $f$-numerically flat if both $E$ and $E^*$ are $f$-nef. Here the transition functions for $E|_S$ are restrictions of those holomorphic transition functions for $E$, and we can canonically identify $E^*|_S$ and $(E|_S)^*$.

\begin{proof}
Let $E$ be $T_X^{1,0}$ or $\Omega_X^{1,0}$. We will prove that $E$ is $f$-nef. Firstly, let $\chi=\hat r\omega_{FS}|_{X_{can}}$ with $\hat r$ a positive rational number, such that $f^*\chi\in2\pi c_1(K_X)$. Then we can uniquely write
$$\omega(t)=e^{-t}\omega_0+(1-e^{-t})f^*\chi+\sqrt{-1}\partial\bar\partial\varphi(t)$$
for a unique family of smooth functions $\varphi(t)$ with $\varphi(0)=0$.

Let $h_t$ be the Hermitian metric on $E$ induced by $\omega(t)$ and $H_t$ the metric on $\mathcal O_{\mathbb P(E)}(1)$ induced by $h_t$. Since $\omega(t)$ is type-III, we can fix a positive number $C$ such that the curvature of $(E,h_t)$
\begin{equation}
\Theta_{h_t}(E)\ge-C\omega(t)\otimes Id_E
\end{equation}
in the sense of Griffiths (see \cite[Chapter III, Section 6]{De}). 
Consider the exact sequence of holomorphic vector bundles over $\mathbb P(E)$:
 $$0\rightarrow F\rightarrow \pi^*E
 \xrightarrow{\hat q}\mathcal O_{\mathbb P(E)}(1)\rightarrow0,$$
 where $\pi:\mathbb P(E)\to X$ is the projection, $\hat q:\pi^*E\to\mathcal O_{\mathbb P(E)}(1)$ is the universal quotient map, and $F:=ker \hat q$, and equip $\pi^*E$ the metric $k_t:=\pi^*h_t$, $\mathcal O_{\mathbb P(E)}(1)$ the induced metric $H_t$ as above, and $F$ the restriction of $k_t$. Up to a unique smooth orthogonal lift $\mathcal O_{\mathbb P(E)}(1)\to\pi^*E$ (which is actually the metric adjoint of $\hat q$), we may regard $F\oplus\mathcal O_{\mathbb P(E)}(1)$ as an orthogonal decomposition of $\pi^*E$ (see \cite[Chapter V, Section 14, p. 273]{De}). Now applying \cite[Chapter V, (14.7), page 275]{De} (also see \cite[Proposition 1.11 (iv)]{DPS}) to the above exact sequence and using the basic fact that $\Theta_{k_t}(\pi^*E)=\pi^*\Theta_{h_t}(E)\ge-C\pi^*\omega(t)\otimes Id_{\pi^*E}$, we conclude 
 $$\Theta_{H_t}(\mathcal O_{\mathbb P(E)}(1))\ge-C\pi^*\omega(t).$$
 
 Next we define a family of Hermitian metrics on $\mathcal O_{\mathbb P(E)}(1)$
 $$\tilde H_t:=H_te^{-C\pi^*\varphi(t)}.$$
 Then it immediately follows that
 \begin{equation}\label{eq-tildeH}
 \Theta_{\tilde H_t}(\mathcal O_{\mathbb P(E)}(1))\ge-C(e^{-t}\pi^*\omega_0+(1-e^{-t})\pi^*f^*\chi).
 \end{equation}
 
Given \eqref{eq-tildeH}, we are now ready to prove $f$-nefness of $E$. To the end, let's fix a reduced irreducible compact analytic subvariety $i:S\hookrightarrow X$ with $f(S)=\{y\}$. Set
$$j:\mathbb P(E|_S)\hookrightarrow \mathbb P(E),\,\,\,\,\pi_S=\pi|_{\mathbb P(E|_S)}:\mathbb P(E|_S)\to S.$$
Naturally, there hold
$$j^*\mathcal O_{\mathbb P(E)}(1)=\mathcal O_{\mathbb P(E|_S)}(1),\,\,\,\,\,\pi\circ j=i\circ\pi_S,$$
and, under a local ambient embedding, the restriction of a smooth metric on $\mathcal O_{\mathbb P(E)}(1)$ is a smooth metric on $\mathcal O_{\mathbb P(E|_S)}(1)$. We claim that 
\begin{equation}\label{eq-tildeH-S}
 \Theta_{\tilde H_t}(\mathcal O_{\mathbb P(E)}(1))|_{\mathbb P(E|_S)}\ge-Ce^{-t}\pi^*\omega_0|_{\mathbb P(E|_S)}
 \end{equation}
Since it is automatically true over $S_{reg}$, let's focus on the claim over singular points of $S$. Let $\iota:X_{can}\hookrightarrow\mathbb{CP}^N$ be the map in the semi-ample fibration \eqref{eq-semiample} and $F:=\iota\circ f$. Choose a coordinate neighborhood $V\subset\mathbb{CP}^N$ of $\iota(y)$, on which there is a smooth real potential $\rho$ with
$$
 \omega_{\mathrm{FS}}|_V=\dd\rho,\qquad \rho(\iota(y))=0.
$$
On the open neighbourhood $U:=F^{-1}(V)$ of $S$, put a smooth function
$$
 u:=\hat r\rho\circ F,\,\,\,\,\, with\,\,\,\,\,f^*\chi|_U=\dd u \,\,\,\,\,and\,\,\,\,\,u|_S=0.
$$
 Define on $\pi^{-1}(U)$ the smooth metric
$$
 \check{H}_t:=\tilde H_te^{-C(1-e^{-t})\pi^*u}.
$$
Because $u|_S=0$, the restriction of $\check{H}_t$ to $\mathbb P(E|_S)$ is exactly the metric occurring in \eqref{eq-tildeH-S}, as
\begin{equation}\label{eq-j-metric}
 j^*\check H_t=j^*\tilde H_t.
\end{equation}
On the ambient open manifold $\pi^{-1}(U)$, the inequality \eqref{eq-tildeH} now gives

 \begin{align}
 \Theta_{\check H_t}(\mathcal O_{\mathbb P(E)}(1))
 &=\Theta_{\tilde H_t}(\mathcal O_{\mathbb P(E)}(1))
     +C(1-e^{-t})\pi^*f^*\chi\nonumber\\
 &\geq-Ce^{-t}\pi^*\omega_0\,\,\,\,\,\,on\,\,\,\,\,U\nonumber.
 \end{align}

This, combining \eqref{eq-j-metric}, concludes that \eqref{eq-tildeH-S} has a smooth ambient representative with the stated lower bound on a neighborhood of every point of $\mathbb P(E|_S)$, including its singular points, and hence \eqref{eq-tildeH-S} follows directly.

Fix a K\"{a}hler metric $\widetilde\omega$ on $\mathbb P(E)$ with $C\pi^*\omega_0\le\widetilde\omega$ on $\mathbb P(E)$. Then \eqref{eq-tildeH-S} implies that the curvature of smooth metrics $j^*\tilde H_t$ on $\mathcal O_{\mathbb P(E|_S)}(1)$ satisfy

\begin{equation}\label{eq-tildeH-S'}
 \Theta_{j^*\tilde H_t}(\mathcal O_{\mathbb P(E|_S)}(1))\ge-Ce^{-t}\pi^*\omega_0|_{\mathbb P(E|_S)}\ge-e^{-t}\widetilde\omega|_{\mathbb P(E|_S)},
 \end{equation}
 which by definition means $\mathcal O_{\mathbb P(E|_S)}(1)$, and hence $E|_S$, is nef. Therefore, $E$ is $f$-nef. Since $E$ can be $T_X^{1,0}$ or its dual $\Omega_X^{1,0}$, we have proved that $T_X^{1,0}$ is $f$-numerically flat.
 
Proposition \ref{prop-f-flat} is proved.
\end{proof}

\subsection{Vanishing of Chern classes}
We shall also show that the existence of a type-III solution on $X$ gives restrictions on the Chern classes $c_k(X)$ of $X$. Precisely, in this subsection we are going to prove the following simple but useful observation.
\begin{prop}\label{prop-chernclass}
Let $X$ be an $n$-dimensional compact connected K\"ahler manifold of a nef $K_X$. Let $\nu(X)$ be the numerical Kodaira dimension of $X$. If there exists a type-III solution $\omega(t)$ to the K\"ahler-Ricci flow \eqref{nkrf} on $X$, then for any $\nu(X)<k\le n$ there holds
$$c_k(X)=0\,\,\,\,in\,\,\,\,H^{2k}(X,\mathbb R).$$
\end{prop}

\begin{proof}[Proof of Proposition \ref{prop-chernclass}]
Recall that the total Chern form with respect to the curvature $\Theta_t$ of the Chern connection on $(T_X^{1,0},\omega(t))$ is defined by
\begin{equation}\label{eq-chernform}
 \det\left(I+\frac{1}{2\pi}\Theta_t\right)
 =\sum_{k=0}^{n}\Gamma_{k,t},
 \,\,\,\,where\,\,\,\, \Gamma_{k,t}=c_k(T_X,\omega(t)).
\end{equation}
Classically we know $\Gamma_{k,t}$ is a real closed $(k,k)$-form representing the Chern class $c_k(T_X^{1,0})$ in real de Rham cohomology. Then it suffices to show that for any smooth $(n-k,n-k)$-form $\xi$ on $X$, 
$$\int_X\Gamma_{k,t}\wedge\xi\to0,\,\,\,\,as\,\,\,t\to\infty.$$
To see this, first note that, since the curvature of $\omega(t)$ is uniformly bounded with respect to $\omega(t)$ and $\Gamma_{k,t}$ is a degree-$k$ homogenerous polynomial in curvature of $\omega(t)$, we may choose a positive number $C_k$ such that
$$|\Gamma_{k,t}|_{\omega(t)}\le C_k, \,\,\,\,\forall t\ge0.$$
Consequently, there is a positive number $C(k,\xi)$, depending only on $k$ and $\xi$, satisfying
$$|\Gamma_{k,t}\wedge\xi|\le C(k,\xi)\omega(t)^k\wedge\omega_0^{n-k}\,\,\,\,on \,\,\,\, X.$$
This imlies, for $k>\nu(X)$,
$$\left|\int_X\Gamma_{k,t}\wedge\xi\right|\le C(k,\xi)\int_X\omega(t)^k\wedge\omega_0^{n-k}=O\left(e^{-(k-\nu(X))t}\right)\to0,\,\,\,\,as\,\,\,t\to\infty.$$
Proposition \ref{prop-chernclass} is proved.
\end{proof}

The above two propositions will be the fundation to prove our main theorems.

\section{Proof of Theorem \ref{thm-proj-a}}\label{sect-proof-proj}
We are ready to give a proof of Theorem \ref{thm-proj-a}. Actually, we are going to prove the following version.

\begin{thm}\label{thm-proj}
Let $X$ be an $n$-dimensional connected projective manifold of a semi-ample canonical line bundle $K_X$ and $1\le m:=kod(X)\le n-1$. Let $f:X\to X_{can}$ be its semi-ample fibration. The followings are equivalent.\\
(1) $X$ admits a type-III solution to the K\"ahler-Ricci flow. \\
(2)There is a commutative diagram
\begin{equation}\label{eq-diagram}
\begin{tikzcd}
 Z \arrow[r,"p"] \arrow[d,"h"'] & X\arrow[d,"f"]\\
 B \arrow[r,"q"'] & X_{can}
\end{tikzcd}
\end{equation}
such that:
\begin{enumerate}
\item [(2.1)] $p$ is a finite surjective \'etale morphism and $q$ is finite and surjective;
\item [(2.2)] $Z$ and $B$ are connected projective manifolds;
\item [(2.3)] $h$ is a holomorphic submersion with connected fibres of dimension $n-m$, and every fiber is a finite \'etale quotient of a complex torus;
\item [(2.4)] $h$ is the semi-ample fibration of $Z$.
\end{enumerate}
\end{thm}

\begin{proof}[Proof of Theorem \ref{thm-proj}, (1) $\Rightarrow$ (2)]
Assume we are in the setup of Theorem \ref{thm-proj}; in particular $X$ is further assumed to be projective. Thanks to Proposition \ref{prop-f-flat}, we know $T_X^{1,0}$ is $f$-nef with respect to the semi-ample fibration $f:X\to X_{can}$. Then applying H\"oring's result \cite[Corollary 2.4]{Ho} gives a finite \'etale cover $p:Z\to X$ for which there is a Stein factorization 
\begin{equation}\label{eq-diagram-1}
\begin{tikzcd}
 Z \arrow[r,"p"] \arrow[d,"h"'] & X\arrow[d,"f"]\\
 B \arrow[r,"q"'] & X_{can},
\end{tikzcd}
\end{equation}
in which $Z$ and $B$ are two normal conneed complex analytic spaces, $h$ is a smooth morphism and $q$ is finite and surjective. Since $X$ is smooth and projective, by the property of a finite \'etale morphism that $p$ is a local biholomorphism between complex analytic spaces, we know $Z$ is also smooth and projective. Moreover, by pulling back an ample line bundle over $X_{can}$ by the finite surjective map $q:B\to X_{can}$, we see that $B$ is also projective. Then the smoothness of $h$ and $Z$ implies that $B$ is also smooth (see The Stacks Project, Lemma 33.25.9, Tag 05AX). Therefore, $h:Z\to B$ is a holomorphic submersion with connected fibers between two connected projective manifolds. 

On the other hand, since $p:Z\to X$ is now a local biholomorphism between two projective manifolds, the differential $dp$ induces a natural isomorphism
\begin{equation}\label{eq-tangentbd}
T_Z^{1,0}=p^*T_X^{1,0}.
\end{equation}
It also directly follows from the commutative diagram \eqref{eq-diagram-1} that for $\tilde S\subset Z$, if $h(\tilde S)=\{pt\}$, then $f(p(\tilde S))=\{pt\}$. Then, combining \eqref{eq-tangentbd}, we conclude that $T_Z^{1,0}$ is $h$-numerically flat. Also $h$ is equidimensional as its fibers are diffeomorphic each other as smooth manifolds. Now we can apply \cite[Lemma 2.1(1)]{Ho} to conclude that each fiber of $h$, which is a smooth fiber, is a finite \'etale quotient of a complex torus. We have proved items (2.1), (2.2) and (2.3) in Theorem \ref{thm-proj}.

Next we prove the final item (2.4) in Theorem \ref{thm-proj}. Fix an integer $b$ and an ample line bundle $\tilde L$ on $X_{can}$ with $bK_X=f^*\tilde L$, and set $L_B=q^*\tilde L$, which is an ample line bundle on $B$ as $q$ is finite and surjective. 

Using \eqref{eq-tangentbd} again, we see
\begin{equation}
bK_Z=p^*(bK_X)=p^*f^*\tilde L=h^*q^*\tilde L=h^*L_B,
\end{equation} 
and hence for any $k\in\mathbb Z_{\ge1}$,
\begin{equation}
K_Z^{\otimes bk}=h^*L_B^{\otimes k}.
\end{equation} 
We claim that 
\begin{equation}\label{eq-h0}
H^0(Z,K_Z^{\otimes bk})=H^0(B,L_B^{\otimes k}).
\end{equation} 
Let's present an elementary proof for this. Firstly recall that $h:Z\to B$ is a surjective holomorphic submersion between smooth projective manifolds, with compact connected fibres. It suffices to prove
$$
 h^*:H^0(B,L)\xrightarrow{\ \sim\ }H^0(Z,h^*L),\,\,\,\,\,with \,\,\,\,\,L:=L_B^k.
$$

\underline{$h^*$ is injective.}
Suppose $s\in H^0(B,L)$ and $h^*s=0$. For any $b\in B$, surjectivity of $h$ gives a point $z\in h^{-1}(b)$. By the definition of a pullback bundle, $(h^*L)_z=L_b$, and in this identification we have $0=(h^*s)(z)=s(b)$. Thus $s=0$. We may mention that this argument also proves injectivity on sections over every open subset of $B$.

\underline{$h^*$ is surjective.}
Let $\tau\in H^0(Z,h^*L)$. Choose an open set $U\subset B$ and a nowhere-zero holomorphic frame $e$ of $L|_U$. The section $h^*e$ is a nowhere-zero frame on the entire inverse image $h^{-1}(U)$, so there is a unique holomorphic function $F$ such that
$$
 \tau|_{h^{-1}(U)}=F\,h^*e,\,\,\,\,\, F\in\mathcal O_Z(h^{-1}(U)).
$$
For each $b\in U$, the restriction $F|_{Z_b}$ is defined on all of the compact connected complex manifold $Z_b=h^{-1}(b)$. Its modulus attains a maximum. In a coordinate ball at a maximum point, the maximum principle makes the function constant, and the identity theorem extends that constant over the connected fiber. Thus $F|_{Z_b}$ is constant.

Define $f:U\to\mathbb C$ by letting $f(b)$ be this constant value. It is well defined, and
$$
 F=f\circ h\,\,\,\,\,on \,\,\,\,\,h^{-1}(U).
$$
We need to check that $f$ is holomorphic. Fix $b_0\in U$ and choose $z_0\in Z_{b_0}$. Set $m=\dim B$, $n=\dim Z$, and choose holomorphic coordinates centred at these points. Since $dh_{z_0}$ is surjective, we may write $h$ as
$$
 h(y,w)=y.
$$
Setting $w=0$ therefore defines a holomorphic map $\sigma:V\to Z$ on a neighbourhood $V\subset U$ of $b_0$, with $\sigma(b_0)=z_0$ and
$$
 h\circ\sigma=Id_V.
$$
Since $\sigma(b)$ lies in the fibre over $b$, the definition of $f$ gives
$$
 f|_V=F\circ\sigma.
$$
This is holomorphic. Such a neighbourhood exists at every $b_0$, hence $f\in\mathcal O_B(U)$. 

Put $s_U=f e\in H^0(U,L)$. Then $h^*s_U=\tau$ over $h^{-1}(U)$. On an overlap of two trivializing open sets $U$ and $V$, the constructed sections satisfy
$$
 h^*(s_U-s_V)=0\,\,\,\,\,on \,\,\,\,\, h^{-1}(U\cap V).
$$
The pointwise injectivity argument above shows $s_U=s_V$ on $U\cap V$. They glue to a unique global holomorphic section $s$ of $L$, and $h^*s=\tau$. This proves surjectivity and completes the proof of \eqref{eq-h0}. 

Next, for sufficiently large $k$ we know $L_B^{\otimes k}$ is very ample. Choose a basis $s_0,\dots,s_N$ of its global sections and let $\iota_k:B\hookrightarrow\mathbb{CP}^N$ be the corresponding embedding. Since these sections have no common zero, their pullbacks have no common zero, too. Thus $K_Z^{\otimes bk}$ is globally generated, proving semi-ampleness of $K_Z$. Moreover, \eqref{eq-h0} says that $h^*s_0,\ldots,h^*s_N$ is a basis of all pluricanonical sections in this degree. In compatible local frames its complete linear-system map is
$$
 z\longmapsto
 [h^*s_0(z),\cdots,h^*s_N(z)]
 =[s_0(h(z)),\cdots,s_N(h(z))]
 =\iota_k(h(z)).
$$
Because $\iota_k$ is an embedding, this map has exactly the fibres of $h$, which are connected. Its image is isomorphic to $B$. Consequently $kod(Z)=dim B$, and $h:Z\to B$ is the semi-ample fibration of $Z$.

The implication (1)$\Rightarrow$(2) in Theorem \ref{thm-proj} is proved.
\end{proof}

\begin{proof}[Proof of Theorem \ref{thm-proj}, (2)$\Rightarrow$(1)]
This direction follows directly from \cite[Theorem 1.5]{ToZyg}, after passing to the finite \'etale cover $Z$. Precisely, for a solution $\omega(t)$ to the K\"{a}hler-Ricci flow \eqref{nkrf} on $X$, set
$$\eta(t):=p^*\omega(t).$$
Then, since $p$ is a local biholomorphism, $\eta(t)$ is a family of K\"{a}hler metrics on the projective manifold $Z$ solving the K\"{a}hler-Ricci flow \eqref{nkrf}, whose curvatures are related by
\begin{equation}\label{eq-rm}
|Rm(\eta(t))|_{\eta(t)}^2(z)=|Rm(\omega(t))|^2_{\omega(t)}(p(z)),\,\,\,\,\,\forall z\in  Z.
\end{equation}
As we have (2), in particular, $Z$ has a semi-ample $K_Z$ and its semi-ample fibration $h:Z\to B$ is a holomorphic submersion between two projective manifolds, whose fiber is a finite \'etale quotient of a complex torus, applying \cite[Theorem 1.5, the last case]{ToZyg} gives that $\eta(t)$ is of type III, i.e. $|Rm(\eta(t))|_{\eta(t)}^2$ is uniformly bounded on $Z\times[0,\infty)$. Using that $p$ is surjective, we conclude from \eqref{eq-rm} that $|Rm(\omega(t))|^2_{\omega(t)}$ is uniformly bounded on $X\times[0,\infty)$, i.e. $\omega(t)$ is of type III.

The implication (2)$\Rightarrow$(1) in Theorem \ref{thm-proj} is proved.
\end{proof}

\begin{rem}\label{rem-torusfibration}
As we have seen in the above implication (2)$\Rightarrow$(1), a submersive semi-ample fibration $h:Z\to B$ with fiber being a finite \'etale quotient of a complex torus suffices to characterize the type-III solutions. Actually, in the setting of Theorem \ref{thm-proj}, one may choose a further commutative diagram
\begin{equation}\label{eq-diagram-torus}
\begin{tikzcd}
 \tilde Z \arrow[r,"r"] \arrow[d,"\tilde h"']
     & Z \arrow[r,"p"] \arrow[d,"h"']
     & X \arrow[d,"f"]\\
 \tilde B \arrow[r,"s"']
     & B \arrow[r,"q"'] & X_{can}
\end{tikzcd}
\end{equation}
in which $r$ and $s$ are finite surjective \'etale morphisms, $\tilde Z$ and $\tilde B$ are connected smooth projective manifolds, and $\tilde h$ is the smooth semi-ample fibration of $\tilde Z$ with connected complex-torus fibres of dimension $n-m$. Setting $\tilde p=p\circ r$ and $\tilde q=q\circ s$ gives a commutative diagram
$$
\begin{tikzcd}
 \tilde Z \arrow[r,"\tilde p"] \arrow[d,"\tilde h"'] & X\arrow[d,"f"]\\
 \tilde B \arrow[r,"\tilde q"'] & X_{can}
\end{tikzcd}
$$
which satisfies all the properties (2.1)-(2.4) in Theorem \ref{thm-proj}, except that the fiber is now strengthened to be a complex torus. Let's prove \eqref{eq-diagram-torus} as follows. Firstly, thanks to H\"oring \cite[Lemma 2.1 (2)]{Ho}, we may choose a finite \'etale cover $r:\tilde Z\to Z$ such that the Stein factorization $\tilde h:\tilde Z\to\tilde B$ is smooth in codimension one and the smooth fibers are tori, and $s:\tilde B\to B$ is finite and surjective. Since $r$ is locally biholomorphic, we know $\tilde Z$ is smooth and projective and $F:=h\circ r:\tilde Z\to B$ is also a holomorphic submersion. Next we shall check that $\tilde B$ is also smooth and $s:\tilde B\to B$ is \'etale by an elmentary argument.

Fix an arbitrary $\tilde b\in\tilde B$ and then set
$$b=s(\tilde b), \,\,\,\,K:=\tilde h^{-1}(\tilde b).$$
For every $x\in K$, since $F:\tilde Z\to B$ is a holomorphic submersion, we may fix a neighborhood $U_x=V_x\times D_x$ of $x$ in $\tilde Z$, where $V_x$ is a neighborhood of $b$ in $B$, $D_x$ is a polydisc in $\mathbb C^{n-m}$ and $x$ corresponds to $(b,0)$, in which $F$ reads
$$F(t,u)=t.$$ 
Note that the connected image $\tilde h(\{t\}\times D_x)\subset s^{-1}(t)$, while $s^{-1}(t)$ is a finite set, it follows that $u\mapsto \tilde h(t,u)$ is constant. Define
$$\tau_x:V_x\to\tilde B,\,\,\,\,\tau_x(t)=\tilde h(t,0),$$
which is a holomorphic map of analytic spaces (as $\tilde h$ is a holomorphic map of analytic spaces by Stein factorization) and satisfies
$$\tilde h(t,u)=\tau_x(t),\,\,\,\,\forall(t,u)\in U_x,\,\,\,\,s\circ\tau_x=Id_{V_x},\,\,\,\,\tau_x(b)=\tilde b.$$

Cover $K$ by finitely many such neighborhoods $U_{x_1},\,...,\,U_{x_N}$ centered at points $x_i\in K$. We have the corresponding base neighborhoods $V_{x_i}\subset B$ and holomorphic maps $\tau_{x_i}:V_{x_i}\to\tilde B$. For $t\in V_{x_i}\cap V_{x_j}$ with $U_{x_i}\cap U_{x_j}\cap K\neq\emptyset$, by their construction we have $\tau_{x_i}(t)=\tau_{x_j}(t)$. Then we may fix a neighborhood $V\subset\cap_{i} V_{x_i}$ and define $\tau:V\to\tilde B$ using $\tau_{x_i}$'s, which satisfies
\begin{equation}
s\circ \tau=Id_{V},\,\,\,\,\,\tau(b)=\tilde b.\nonumber
\end{equation}
Define
$$W:=\bigcup(U_{x_i}\cap F^{-1}(V)).$$
This is a neighborhood of $K$, contained in $F^{-1}(V)$, and
$$\tilde h|_W=\tau\circ F|_W.$$
It easily follows that on the neighborhood $N:=\tilde B\setminus\tilde h(\tilde Z\setminus W)$ of $\tilde b\in\tilde B$,
$$y=\tau(s(y)),\,\,\,\forall y\in N,$$
implying that
$$N\subset s^{-1}(V),\,\,\,\,\,\tau\circ s|_N=Id_N.$$
Set $V_0:=\tau^{-1}(N)\subset V$, then it satisfies $s(N)=V_0$ and
$$s|_N:N\to V_0,\,\,\,\,\,\tau|_{V_0}:V_0\to N$$
are holomorphic inverses. In particular, since $\tilde b\in \tilde B$ is chosen arbitrarily and $B$ is smooth, we see that $\tilde B$ is a smooth complex manifold and $s:\tilde B\to B$ is a local biholomorphism and hence a finite \'etale cover. This also implies that $\tilde B$ is projective by pulling back an ample line bundle on $B$ to $\tilde B$ by $s$. Finally, recall that $h:Z\to B$ is a holomorphic submersion, we conclude that $\tilde h:\tilde Z\to\tilde B$ is a holomorphic submersion, and hence each fiber is smooth and is an $(n-m)$-dimension complex torus. Also, we may identify $\tilde h$ with the semi-ample fibration of $\tilde Z$ by the same argument in the last part of the proof of Theorem \ref{thm-proj}, (1)$\Rightarrow$(2). We have proved all the claimed properties for diagram \eqref{eq-diagram-torus}.
\end{rem}

\section{Proof of Theorem \ref{thm-kahler-a}}\label{sect-proof-kahler}
We now prove the following more precise version of Theorem \ref{thm-kahler-a}.

\begin{thm}\label{thm-kahler}
Let $X$ be an $n$-dimensional compact connected K\"ahler manifold of a nef canonical line bundle $K_X$ and the numerical Kodaira dimension $\nu(X)=1$. The followings are equivalent.\\
(1) $X$ admits a type-III solution to the K\"ahler-Ricci flow. \\
(2) $K_X$ is semi-ample, and its semi-ample fibration $f:X\to X_{can}$ satisfies a commutative diagram
\begin{equation}\label{eq-diagram-kod1}
\begin{tikzcd}
 Z \arrow[r,"p"] \arrow[d,"h"'] & X\arrow[d,"f"]\\
 B \arrow[r,"q"'] & X_{can}
\end{tikzcd}
\end{equation}
such that:
\begin{enumerate}
\item [(2.1)] $p$ is a finite surjective \'etale morphism and $q$ is finite and surjective;
\item [(2.2)] $Z$ is a compact connected K\"ahler  manifold and $B$ is a closed Riemann surface of genus at least two;
\item [(2.3)] $h$ is a holomorphic submersion with connected fibres of dimension $n-1$, and every fiber is an $(n-1)$-dimensional complex torus;
\item [(2.4)] $h$ is the semi-ample fibration of $Z$.
\end{enumerate}
\end{thm}

\begin{proof}[Proof of Theorem \ref{thm-kahler}]
As in Theorem \ref{thm-proj}, the implication (2)$\Rightarrow$(1) follows directly from \cite[Theorem 1.5]{ToZyg} (also see \cite{FZ,HeTo}). 

To see the implication (1)$\Rightarrow$(2), we firstly apply Proposition \ref{prop-chernclass} and the assumption $\nu(X)=1$ to see that 
$c_2(X)=0\,\,\,\,in\,\,\,\,H^{4}(X,\mathbb R).$
In particular, $c_2(X)=0$ in $H^{2,2}(X,\mathbb R)$. Now the abundance theorem of Iwai-Matsumura \cite[Theorems 1.2 and 1.3]{IM} can be applied to conclude that $\Omega_X^{1,0}$ is nef and $K_X$ is semi-ample; then we further apply \cite[Theorem 3.5]{IM} to obtain a finite \'etale cover $p:Z\to X$ of $X$ with $K_Z$ semi-ample, and the semi-ample fibration of $Z$,
$$h:Z\to B,$$
is a holomorphic submersion over a Riemann surface of genus at least two, whose fibers are complex tori, and there is a finite surjection $q:B\to X_{can}$ such that these maps satisfy the conclusions in item (2) in Theorem \ref{thm-kahler}. This finishes the proof of the implication (1)$\Rightarrow$(2).

Theorem \ref{thm-kahler} is proved.
\end{proof}

\begin{rem}\label{rem-onekod}
Basing the commutative diagram \eqref{eq-diagram-kod1} and that $X_{can}$ and $B$ are closed Riemann surfaces, we now give a more explicit description for the singular fibers of the semi-ample fibration $f:X\to X_{can}$ (compare \cite[Proof of Lemma 2.1, step 3]{Ho} and \cite[Theorem 1.6]{ToZyg}). Fix a $s\in V\subset X_{can}$ be a critical value of $f$, $x\in f^{-1}(s)$ and $x'\in p^{-1}(x)$, and then set $s'=h(x')$. Since $q:B\to X_{can}$ is a holomorphic map between two Riemann surface, by one-variable complex analysis theory we may choose coordinates $v$ near $s$ and $u$ near $s'$, centered at $s$ and $s'$ respectively, such that
$$v\circ q(u)=u^e, \,\,\,\,e\ge1.$$
Also, since $h:Z\to B$ is a holomorphic submersion, we may further fix a coordinate $(z^1,\,...,\,z^n)$ centered at $x'$ such that
$$u\circ h(z)=z^1.$$
Then applying the commutative diagram \eqref{eq-diagram-kod1} gives
$$v\circ f\circ p(z)=v\circ q\circ h(z)=(u\circ h(z))^e=(z^1)^e.$$
Since $p$ is local biholomorphic, one may fix neighborhoods $x\in U\subset X$ and $x'\in U'\subset Z$ such that $p|_{U'}:U'\to U$ is biholomorphic, and then define coordinates $w=z\circ(p|_{U'})^{-1}$ on $U$, which satisfies
$$v\circ f(w)=(w^1)^e\,\,\,\,\,on\,\,\,\,\, U.$$
This means that locally near an arbitrary $x\in f^{-1}(s)$, as a subset in $X$, $f^{-1}(s)$ is a coordinate hyperplane $\{w^1=0\}$, while as a divisor, $f^{-1}(s)$ is a coordinate hyperplane with coefficient $e$. Moreover, when $x\in f^{-1}(s)$ varies, we know $e=e(x)$ is locally constant and hence constant on $f^{-1}(s)$, as $f^{-1}(s)$ is connected. In conlusion, the singular fiber $X_s:=f^*s$ over $s$ is given by
$$X_s=m_sF_s,$$
where $F_s$ is an $(n-1)$-dimensional connected smooth complex manifold and $m_s=e$, which must be at least $2$ as $s$ is a critical value of $f$. Moreover, given the above construction, we have that, when restricting to the fiber $Z_{s'}$ of $h$ over $s'$, $p|_{Z_{s'}}:Z_{s'}\to F_s$ is locally given by $(z^2,\,...,\,z^n)\mapsto(w^2,\,...,\,w^n)$, and hecne it is a finite local biholomorphism, namely, $F_s$ is a finite \'etale quotient of the complex torus $Z_{s'}$. 

Finally, as an application, we mention that $Id_X:X\to X$ is already a log-resolution of $X_s$ in $X$, and hence its log-canonical threshold is $\frac{1}{m_s}$ and its log-canonical multiplicity is $1$ (see their defintions in \cite[Section 2]{Zys20b}), and hence we have recovered \cite[Theorem 1.3]{Zys20b}.
\end{rem}

\section{Proof of Theorem \ref{thm-kahler-three-a}}\label{sect-proof-kahler-three}
We now prove the following more precise version of Theorem \ref{thm-kahler-three-a}.

\begin{thm}\label{thm-kahler-three}
Let $X$ be a $3$-dimensional compact connected K\"ahler manifold of a semi-ample canonical line bundle $K_X$ and $kod(X)=2$. The followings are equivalent.\\
(1) $X$ admits a type-III solution to the K\"ahler-Ricci flow. \\
(2) Its semi-ample fibration $f:X\to X_{can}$ satisfies a commutative diagram
\begin{equation}\label{eq-diagram-three}
\begin{tikzcd}
 Z \arrow[r,"p"] \arrow[d,"h"'] & X\arrow[d,"f"]\\
 B \arrow[r,"q"'] & X_{can}
\end{tikzcd}
\end{equation}
such that:
\begin{enumerate}
\item [(2.1)] $p$ is a finite surjective \'etale morphism and $q$ is finite and surjective;
\item [(2.2)] $Z$ is a $3$-dimensional compact connected K\"ahler manifold and $B$ is a smooth connected projective surface;
\item [(2.3)] $h$ is a holomorphic submersion with connected fibres and every fiber is an $1$-dimensional complex torus;
\item [(2.4)] $h$ is the semi-ample fibration of $Z$.
\end{enumerate}
\end{thm}

\begin{proof}[Proof of Theorem \ref{thm-kahler-three}]
It suffices to prove the implication (1)$\Rightarrow$(2).

By Proposition \ref{prop-f-flat} (or more directly by the semi-ample fibration, $K_X$ is $f$-numerically flat. Recall \cite[Proposition 1.4]{ToZyg} that each fiber of $f$ contains no rational curves, and hence we apply \cite[Lemma 3.1]{Ho} to see that $f$ does not contract any divisor, namely $f$ is equidimensional. Using again Proposition \ref{prop-f-flat} and \cite[Lemma 2.1(2)]{Ho}, there is a finite \'etale cover $p_0:Z_0\to X$ such that the Stein factorization $h_0:Z_0\to B_0$ is smooth in codimension one and the smooth fibers are one-dimensional tori, and we have a commutative diagram
\begin{equation}\label{eq-diagram-three-1}
\begin{tikzcd}
 Z_0 \arrow[r,"p_0"] \arrow[d,"h_0"'] & X\arrow[d,"f"]\\
 B_0 \arrow[r,"q_0"'] & X_{can}
\end{tikzcd}
\end{equation}
where $q_0$ is finite and surjective, $Z_0$ is a $3$-dimensional compact connected K\"ahler manifold (since $X$ is smooth K\"ahler and $p_0$ is a finite \'etale cover), $B_0$ is a normal projective surface (as pulling back an ample line bundle on $X_{can}$ by the finite surjective map $q_0$ gives an ample line bundle on $B_0$). Moreover, since $h_0$ is smooth in codimension one and $B_0$ is normal, we  may fix a finite subset $S_0\subset B_0$, such that $B_0\setminus S_0$ is smooth and $h_0:Z_0\setminus h_0^{-1}(S_0)\to B_0\setminus S_0$ is a holomorphic submersion with fibers one-dimensional tori (i.e. elliptic curves). Since $S_0$ has codimension two in $B_0$, it is a basic fact (see \cite[Lemma 3 and its proof]{An} or \cite[Proof of Theorem 1.4]{Ho}) that $h_0:Z_0\setminus h_0^{-1}(S_0)\to B_0\setminus S_0$ is isotrivial, namely, $h_0:Z_0\to B_0$ has generically constant moduli. On the other hand, because $p_0$ is a local biholomorphism, we easily see that $h_0$ is also equidimensional and $T_{Z_0}$ is $h_0$-numerically flat. Then we apply \cite[Lemma 2.1 (4)]{Ho} to $h_0:Z_0\to B_0$ to conclude a finite \'etale cover $p_1:Z\to Z_0$ with the smooth and locally trivial Stein factorization $h:Z\to B$, and there is a commutative diagram
\begin{equation}\label{eq-diagram-three-2}
\begin{tikzcd}
 Z \arrow[r,"p_1"] \arrow[d,"h"'] & Z_0\arrow[d,"h_0"]\\
 B \arrow[r,"q_1"'] & B_0
\end{tikzcd}
\end{equation}
where $q_1$ is finite and surjective, $Z$ is a $3$-dimensional compact connected K\"ahler manifold, $B$ is a smooth projective surface. In particular, $h:Z\to B$ is a holomorphic submersion, and it is the semi-ample fibration of $Z$, using the same arguments in the last part of proof of Theorem \ref{thm-proj}, (1)$\Rightarrow$(2). Setting $p:=p_0\circ p_1:Z\to X$ and $q:=q_0\circ q_1$, they satisfy the commutative diagram \eqref{eq-diagram-three} and we have conmpleted the proof of Theorem \ref{thm-kahler-three}, (1)$\Rightarrow$(2).
\end{proof}

\appendix

\section{Examples of type-III solution to the K\"{a}hler-Ricci flow on manifolds of nonpositive holomorphic curvature}
In this appendix, we exhibit examples of type-III solution to the K\"{a}hler-Ricci flow, which indicates that, for a compact K\"{a}hler manifold, the existence of a type-III solution to the K\"{a}hler-Ricci flow may be regarded as a curvature condition generalizing the existence of a K\"{a}hler metric with nonpositive holomorphic curvature.

\subsection{Case of nonpositive bisectional curvature} 
Let $X$ be an $n$-dimensional compact connected K\"{a}hler manifold admitting a K\"{a}hler metric with nonpositive bisectional curvature. We shall construct a type-III solution to the K\"{a}hler-Ricci flow on $X$, using the structure theorem of such manifold. Precisely, thanks to Liu's work \cite[Corollary 4]{Liu}, there is a canonical K\"ahler metric $\omega_*$ on $X$ locally equal to $\omega_F+\omega_E$, where $\omega_F$ is flat of dimension $n-m$, and $Ric(\omega_E)=-\omega_E$ on the other $m$ dimensions ($m=kod(X)$).

Define \emph{globally} on $X$
\begin{equation}\label{eq-ricform}
 \theta=-Ric(\omega_*),\,\,\,\,\, \eta=\omega_*+Ric(\omega_*).
\end{equation}
In the local product charts where $\omega_*=\omega_F+\omega_E$, these are precisely $\theta=\omega_E$ and $\eta=\omega_F$. Thus they are closed semipositive real $(1,1)$-forms and $\eta+\theta=\omega_*>0$. Their intrinsic definitions in \eqref{eq-ricform} show that they agree on overlapping charts, in particular, they are globally well-defined on $X$.

For $t\ge0$, set
\begin{equation}\label{eq-flow1}
 \omega_*(t)=e^{-t}\eta+\theta.
\end{equation}
This is positive for every finite $t$. In local product coordinates its matrix is block diagonal, with blocks $e^{-t}\omega_F$ and $\omega_E$, so
$$
\omega_*(t)^n=\binom{n}{m}e^{-(n-m)t} \omega_F^{n-m}\wedge\omega_E^m.
$$
and hence
$$
 Ric(\omega_*(t))=Ric(\omega_E)=-\theta,\,\,\,\,\,
 \partial_t\omega_*(t)=-e^{-t}\eta
   =-Ric(\omega_*(t))-\omega_*(t).
$$
Thus \eqref{eq-flow1} solves the K\"{a}hler-Ricci flow \eqref{nkrf} on $X$ with initial metric $\omega_*$.

The full curvature bound follows from the same local product structure. There are no mixed curvature components, and all components in the shrinking flat factor vanish. The only nonzero components are those of the fixed $\omega_E$ factor, whose metric also remains fixed. Therefore,
\begin{equation}
 \sup_{X\times[0,\infty)}|Rm(\omega_*(t))|_{\omega_*(t)}=\sup_X|Rm(\omega_*)|_{\omega_*}<\infty.\nonumber
\end{equation}
This explicitly checks that $\omega_*(t)$ is a type-III solution on $X$.

\subsection{Case of nonpositive holomorphic sectional curvature} We now switch to the case of nonpositive holomorphic sectional curvature.

\subsubsection{Projective case} Let's first look at the projective case. Here we introduce a notation: for a projective manifold $M$ of nef $K_M$, write $n(M)$ for the nef dimension of $M$ defined in \cite[Theorem 2.4]{HLW16}, which may be different from the numerical Kodaira dimension $\nu(M)$. If $K_M$ is semi-ample, then there holds $n(M)=\nu(M)=kod(M)$.

Let $X$ be an $n$-dimensional connected \emph{projective} manifold of a semi-ample canonical line bundle, which admits a K\"{a}hler metric with nonpositive holomorphic sectional curvature. Heier-Lu-Wong-Zheng's splitting theorem \cite[Theorem 1.6]{HLWZ} gives a finite \'etale cover
\begin{equation}
 p:Z=F\times B\longrightarrow X, \,\,\,\,\,dim B=n(B)=n(X),\nonumber
\end{equation}
where $F$ is a projective complex torus and $B$ is a smooth projective manifold admitting a K\"ahler metric of nonpositive holomorphic sectional curvature. May assume $dim B\ge1$.	Denote the factor projections by $\pi_F$ and $\pi_B$.

Since $p$ is a finite \'etale cover and $K_F$ is a trivial line bundle,
$$
 K_Z\simeq p^*K_X\simeq\pi_B^*K_B.
$$
Recall that we have assumed that $K_X$ is semi-ample, and note that $p$ is locally biholomorphic, we know $K_Z$ is also semi-ample. Fix an integer $k\ge1$ such that $K_Z^{\otimes k}$ is base-point-free, then restricting to $\{a\}\times B$  shows that $K_B^{\otimes k}$ is base-point-free, too; namely, $K_B$ is semi-ample. Then applying \cite[Corollary 1.8]{HLWZ} gives
$$
 kod(B)=\nu(B)=n(B)=dim B.
$$
Consequently, the manifold $B$ is a minimal model of general type. Further, the existence of a K\"ahler  metric with nonpositive holomorphic sectional curvature excludes rational curves \cite{Yau}, then by a classical result of Kawamata \cite{Ka} that a minimal model of general type without rational curves has ample canonical line bundle, we conclude that $K_B$ is ample, i.e.
$$
c_1(B)<0.
$$
Now for any solution $\omega(t)$ to the K\"ahler-Ricci flow \eqref{nkrf} on $X$, $p^*\omega(t)$ is a solution to the K\"ahler-Ricci flow \eqref{nkrf} on $Z$, which has uniformly bounded curvature by Gill \cite{Gill}. As before, this implies that $\omega(t)$ has uniformly bounded curvature, namely, $\omega(t)$ is a type-III solution.

\subsubsection{Compact K\"ahler case} Next, we consider the compact K\"ahler case. Let $X$ be an $n$-dimensional compact connected \emph{K\"ahler} manifold of a semi-ample canonical line bundle and $1\le m:=kod(X)\le n-1$, which admits a K\"{a}hler metric with nonpositive holomorphic sectional curvature. For such $X$, we can't apply \cite[Theorem 1.6]{HLWZ} to conclude a global product finite \'etale cover for $X$. However, by following \cite[Proof of Theorem 1.6, page 957-959]{HLWZ} (also see \cite{WuZh}) step by step, and replacing the role of nef reduction map therein by semi-ample fibration map $f:X\to X_{can}$, we can still obtain a finite \'etale cover $p:Z\to X$ and a holomorphic fiber bundle $h:Z\to B$ (which is exactly the map $h_1:M'\to B$ in \cite[page 959]{HLWZ}), where $Z$ and $B$ are compact connected K\"{a}hler manifolds of dimensions $n$ and $m$, respectively, $B$ admits a K\"{a}hler metric of nonpositive holomorphic sectional curvature, and fiber of $h$ is biholomorphic to a fixed $(n-m)$-dimensional complex torus.

To produce a type-III solution on $X$, it suffices to further check that $h:Z\to B$ is the semi-ample fibration of $Z$. Firstly, because $p:Z\to X$ is a finite \'etale cover, which preserves Kodaira dimension \cite[Theorem 6.10]{U75}, and $K_X$ is sem-ample, we see that $K_Z$ is semi-ample and $kod (Z)=kod(X)=m$. Next, we recall from the construction in \cite[page 959]{HLWZ} that $Z$ has a universal cover $\mathbb C^{n-m}\times B'$ and each deck transformation acts diagonally, the $(1,0)$-form frame $dz^1,\,...,\,dz^{n-m}$ on $\mathbb C^{n-m}$ decends to global nowhere-vanishing holomorphic $(1,0)$-forms $\lambda^1$, ..., $\lambda^{n-m}$ on $Z$. Locally over a small open neighborhood $U\subset B$, we know $h^{-1}(U)=T\times B$ holomorphically, and if $\beta$ is a local frame of $K_B$, then 
$$\lambda^1\wedge\cdots\wedge\lambda^{n-m}\wedge h^*\beta$$
is a local frame of $K_Z$ over $h^{-1}(U)$, which transforms exactly as $\beta$. Therefore, we conclude
$$K_Z=h^*K_B, \,\,\,\,\,H^0(Z,K_Z^{\otimes k})=H^0(B,K_B^{\otimes k}).$$
Consequently, $K_B$ is semi-ample and $kod(B)=kod(Z) =m=dim B$. Now $B$ is a compact K\"ahler manifold of a nef and big canonical line bundle $K_B$, then by Moishezon's theorem it is a projective manifold. Further, $B$ contains no rational curves, and hence $K_B$ is ample. Then we can apply the same discussion in the last part of proof of Theorem \ref{thm-proj}, (1)$\Rightarrow$(2), to conclude that $h:Z\to B$ is the semi-ample fibration of $Z$. Given the above structure, we can apply the same discussion in proof of Theorem \ref{thm-proj}, (2)$\Rightarrow$(1), to see that every solution to the K\"ahler-Ricci flow on $X$ is type-III.

\begin{rem}
In general, even in the case of projective manifolds of a semi-ample canonical line bundle, admitting a type-III solution to the K\"ahler-Ricci flow is a curvature condition \emph{strictly} weaker than admitting a K\"ahler metric of nonpositive holomorphic sectional curvature. In two dimension, a simple example is provided by the Fermat quintic surface $
B=\{z_0^5+z_1^5+z_2^5+z_3^5=0\}\subset\mathbb P^3$, which has ample $K_B$ and contains rational curves. On the one hand, given the ampleness of $K_B$, by Cao's convergence theorem \cite{C} we know there is a type-III solution to the K\"ahler-Ricci flow on $B$. On the other hand, the existence of rational curves on $B$ rules out the existence of a K\"ahler metric of nonpositive holomorphic sectional curvature. For higher dimensional examples, consider a product K\"ahler manifold $X=F\times B$, where $F$ is a complex torus. Then $X$ is a compact K\"ahler manifold (which is projective when $F$ is projective) of a semi-ample canonical line bundle, and it admits a type-III solution to the K\"ahler-Ricci flow (by \cite{C,Gill}) but no K\"ahler metric with nonpositive holomorphic sectional curvatur as $X$ contains rational curves.
\end{rem}

\section*{AI Disclosure}
LLMs are used as auxiliary tools to check some of the mathematical computations and arguments. The author takes full responsibility for the content of this paper, which is independently written and verified by the author and contains no AI-generated content.

\end{document}